\documentclass[10pt]{amsart}

\usepackage{amsmath,amsfonts,amssymb,amsthm,epsfig,young}
\usepackage[vcentermath]{youngtab}
\usepackage{hyperref}

\def\C{\mathbb{C}}
\def\Z{\mathbb{Z}}
\def\N{\mathbf{N}}

\def\g{\ensuremath{\mathfrak{g}}}

\def\V{\mathbf{V}}

\def\v{\mathbf{v}}
\def\w{\mathbf{w}}

\def\ke{{\tilde e}}
\def\kf{{\tilde f}}

\DeclareMathOperator{\im}{Im} 
\DeclareMathOperator{\Hom}{Hom}

\DeclareMathOperator{\End}{End}

\DeclareMathOperator{\inc}{in}
\DeclareMathOperator{\out}{out}
\DeclareMathOperator{\tr}{tr}

\DeclareMathOperator{\wt}{wt}

\DeclareMathOperator{\Coker}{Coker}

\newtheorem{theo}{Theorem}[section]
\newtheorem{prop}[theo]{Proposition}
\newtheorem{lem}[theo]{Lemma}
\newtheorem{cor}[theo]{Corollary}

\newtheorem{defin}[theo]{Definition}
\newtheorem*{rem*}{Remark}

\numberwithin{equation}{section}

\begin{document}
\title{Geometric and Combinatorial Realizations of Crystals of Enveloping Algebras}
\author{Alistair Savage}
\address{Fields Institute and University of Toronto \\
Toronto, Ontario \\ Canada}
\email{alistair.savage@aya.yale.edu}
\thanks{This research was supported in part by the Natural
Sciences and Engineering Research Council (NSERC) of Canada}
\date{December 30, 2005}

\begin{abstract}
Kashiwara and Saito have defined a crystal structure on the set of
irreducible components of Lusztig's quiver varieties.  This gives a
geometric realization of the crystal graph of the lower half of the
quantum group associated to a simply-laced Kac-Moody algebra.  Using
an enumeration of the irreducible components of Lusztig's quiver
varieties in finite and affine type $A$ by combinatorial data, we
compute the geometrically defined crystal structure in terms of this
combinatorics.  We conclude by comparing the combinatorial
realization of the crystal graph thus obtained with other
combinatorial models involving Young tableaux and Young walls.
\end{abstract}

\maketitle

\section*{Introduction}
Crystal graphs may be viewed as the $q \to 0$ limit of the theory of
quantum groups.  In this limit, representation theory becomes
combinatorics.  A representation is replaced by its crystal graph,
from which important information such as its character can be
computed.  Other information such as the decomposition of tensor
products of representations into irreducible ones can also be
determined from the crystal graphs of the representations involved.
It is thus an important problem to have concrete realizations of
crystal graphs.

There are many purely combinatorial realizations of crystal graphs
(see \cite{Kas03}). These include (among others) constructions
involving Young tableaux and Young walls, Littelmann's path model,
and the Kyoto path realization.  In \cite{KS97}, Kashiwara and Saito
described a geometric realization of the crystals of the lower half
$U_q^-(\g)$ of the quantum group corresponding to a simply-laced
Kac-Moody algebra $\g$ (a $q$-deformation of the lower half of the
enveloping algebra of $\g$). In \cite{S02}, Saito gave a similar
realization of the crystals of the irreducible integrable
highest-weight representations of these quantum groups. These
constructions involve the quiver varieties of Lusztig and Nakajima.
More specifically, the underlying set of the crystal is the set of
irreducible components of certain varieties associated to quivers
and the operators defining the crystal structure are described in a
natural geometric way.

In \cite{Sav03b}, the author described an explicit isomorphism
between the geometric realization of crystal graphs of integrable
representations and the purely combinatorial approaches involving
Young tableaux and Young walls.  In the current paper, we turn our
attention to the crystals of the lower half of the quantum groups.
Enumerating the irreducible components of the quiver varieties of
finite and affine type $A$ by combinatorial data, we describe the
geometric crystal action in terms of this data.  We also describe
the relationship between the combinatorial model thus obtained and
recently developed models involving Young tableaux and Young walls.

The geometric approach to crystals has several advantages over
purely combinatorial ones. First of all, the quiver variety
construction works uniformly for all simply-laced Kac-Moody algebras
whereas many combinatorial models are often particular to certain
subclasses. For instance, the Young tableaux approach works for the
semi-simple Lie algebras and the Young wall approach for the affine
Lie algebras. Secondly, the geometric viewpoint often yields
alternative (frequently simpler) direct proofs of the combinatorial
constructions.  Once the irreducible components are enumerated by
the desired combinatorial data, one simply computes the action of
the crystal operators.  We know in advance what crystal we have and
thus avoid the need to prove that the combinatorial model yields the
crystal we want. Finally, and perhaps most importantly, the theory
of quiver varieties allows one to ``thaw'' the crystals.  That is,
there is a general geometric procedure for constructing the full
representations and not just the $q \to 0$ limit.  In this
construction, a basis of the representation in question is naturally
enumerated by the irreducible components of the quiver variety. From
purely combinatorial constructions, it is often a highly non-trivial
task to recover the full representation.  However, if we know the
relationship between the geometric and combinatorial realizations,
we can use the general theory of quiver varieties to thaw the
combinatorial crystal.

We expect that the results of this paper can be extended to types
other than $A$.  For instance, in \cite{Sav03b}, types $A$ and $D$
were considered for the crystals of highest-weight integrable
representations. However, since the proofs in type $A$ are simpler,
we restrict ourselves to this case in the current paper for
expository reasons.

The organization of this paper is as follows. In
Sections~\ref{sec:lus_def} and~\ref{sec:ca_qv} we recall the
definition of Lusztig's quiver varieties and the crystal structure
on their irreducible components.  We give a combinatorial
enumeration of these components in Section~\ref{sec:enumerate-comps}
and compute the geometric crystal action in terms of this
combinatorics in Section~\ref{sec:comb-crystals}.  Finally, in
Section~\ref{sec:comparison}, we compare our combinatorial model to
the ones in terms of Young tableaux and Young walls.


\section{Lusztig's quiver variety}
\label{sec:lus_def}

In this section, we will recount the description given in
\cite{L91} of Lusztig's quiver variety and its irreducible
components. See this reference for details, including proofs.

Let $I$ be the set of vertices of the Dynkin graph of a symmetric
Kac-Moody Lie algebra $\mathfrak{g}$ and let $H$ be the set of
pairs consisting of an edge together with an orientation of it.
For $h \in H$, let $\inc(h)$ (resp. $\out(h)$) be the incoming
(resp. outgoing) vertex of $h$.  We define the involution $\bar{\
}: H \to H$ to be the function which takes $h \in H$ to the
element of $H$ consisting of the same edge with opposite
orientation.  An \emph{orientation} of our graph is a choice of a
subset $\Omega \subset H$ such that $\Omega \cup \bar{\Omega} = H$
and $\Omega \cap \bar{\Omega} = \emptyset$.

Let $\mathcal{V}$ be the category of
finite-dimensional $I$-graded vector spaces $\V = \oplus_{i
  \in I} \V_i$ over $\C$ with morphisms being linear maps
respecting the grading.  Then $\V \in \mathcal{V}$ shall denote that
$\V$ is an object of $\mathcal{V}$.

Given $\V \in \mathcal{V}$, let
\[
\mathbf{E_V} = \bigoplus_{h \in H} \Hom (\V_{\out(h)},
\V_{\inc(h)}).
\]
For any subset $H' \subset H$, let $\mathbf{E}_{\V, H'}$ be the
subspace of $\mathbf{E_V}$ consisting of all vectors $x = (x_h)$
such that $x_h=0$ whenever $h \not\in H'$.  The algebraic group
$G_\V = \prod_i GL(\V_i)$ acts on $\mathbf{E_V}$ and
$\mathbf{E}_{\V, H'}$ by
\[
(g,x) = ((g_i), (x_h)) \mapsto (g_{\inc(h)} x_h g_{\out(h)}^{-1}).
\]

Define the function $\varepsilon : H \to \{-1,1\}$ by $\varepsilon
(h) = 1$ for all $h \in \Omega$ and $\varepsilon(h) = -1$ for all $h
\in {\bar{\Omega}}$.  Let $\left<\cdot,\cdot\right>$ be the
non-degenerate, $G_\V$-invariant, symplectic form on $\mathbf{E_V}$
with values in $\C$ defined by
\[
\left<x,y\right> = \sum_{h \in H} \varepsilon(h) \tr (x_h y_{\bar{h}}).
\]
Note that $\mathbf{E_V}$ can be considered as the cotangent space of
$\mathbf{E}_{\V, \Omega}$ under this form.

The moment map associated to the $G_{\V}$-action on the symplectic
vector space $\mathbf{E_V}$ is the map $\psi : \mathbf{E_V} \to
\mathbf{gl_V} = \prod_i \End(\V_i)$ with $i$-component $\psi_i :
\mathbf{E_V} \to \End \V_i$ given by
\[
\psi_i(x) = \sum_{h \in H,\, \inc(h)=i} \varepsilon(h) x_h x_{\bar{h}} .
\]

\begin{defin}[\cite{L91}]
\label{def:nilpotent} An element $x \in \mathbf{E_V}$ is said to be
\emph{nilpotent} if there exists an $N \ge 1$ such that for any
sequence $h_1, h_2, \dots, h_N$ in $H$ satisfying $\out (h_1) = \inc
(h_2)$, $\out (h_2) = \inc (h_3)$, \dots, $\out (h_{N-1}) = \inc
(h_N)$, the composition $x_{h_1} x_{h_2} \dots x_{h_N} : \V_{\out
(h_N)} \to \V_{\inc (h_1)}$ is zero.
\end{defin}

\begin{defin}[\cite{L91}] Let $\Lambda_\V$ be the set of all
  nilpotent elements $x \in \mathbf{E_V}$ such that $\psi_i(x) = 0$
  for all $i \in I$.
\end{defin}

We call $\Lambda_\V$ \emph{Lusztig's quiver variety}. Since
$\Lambda_\V$ depends, up to isomorphism, only on the graded
dimension of $\V$, we define
\[
\Lambda(\v) = \Lambda_{\V^\v},\quad \text{where }(\V^\v)_i =
\C^{\v_i}.
\]


\section{Crystal action on quiver varieties}
\label{sec:ca_qv}

In this section, we review the realization of the crystal graph of
integrable highest-weight representations of a Kac-Moody algebra \g\
with symmetric Cartan matrix via quiver varieties. For details, we
refer the reader to \cite{KS97}.

Let $\mathbf{v, \bar \v, v'} \in (\Z_{\ge 0})^I$ be such that $\v =
\bar \v + \v'$.  Consider the maps
\begin{equation}
\label{eq:diag_action} \Lambda(\v') \times \Lambda(\bar \v)
\stackrel{p_1}{\leftarrow} \tilde \Lambda (\v',\bar \v)
\stackrel{p_2}{\rightarrow} \Lambda(\v).
\end{equation}
Here $\tilde \Lambda (\v',\bar \v)$ is the variety of tuples
$(x,\phi',\bar \phi)$ where $x \in \Lambda(\v)$ and $\phi' =
(\phi'_i)_{i \in I}$, $\bar \phi = (\bar \phi_i)_{i \in I}$ give an
exact sequence
\[
0 \to V^{\v'} \stackrel{\phi'}{\longrightarrow} V^\v \stackrel{\bar
\phi}{\longrightarrow} V^{\bar \v} \to 0
\]
such that $\im \phi'$ is $x$-stable, that is, $x_h ((\im
\phi')_{\out(h)}) \subseteq (\im \phi')_{\inc(h)}$ for all $h \in
H$.  Then $x$ induces $x' \in \Lambda(\v')$ and $\bar x \in
\Lambda(\bar \v)$. Note that we have used the fact that $x$ is
nilpotent if and only if the induced $x'$ and $\bar x$ are. The maps
in \eqref{eq:diag_action} are defined by $p_1(x,\phi',\bar \phi) =
(x',\bar x)$, $p_2(x,\phi',\bar \phi) = x$.

For $i \in I$, define $\varepsilon_i : \Lambda(\v) \to \Z_{\ge 0}$
by
\[
\varepsilon_i(x) = \dim_\C \Coker \left( \bigoplus_{h\, :\,
\inc(h)=i} V_{\out(h)} \stackrel{(x_h)}{\longrightarrow} V_i
\right).
\]
Then, for $c \in \Z_{\ge 0}$, let
\[
\Lambda(\v)_{i,c} = \{x \in \Lambda(\v)\ |\ \varepsilon_i(x) = c\}.
\]
It is easily seen that $\Lambda(\v)_{i,c}$ is a locally closed
subvariety of $\Lambda(\v)$.

From now on, let $\bar \v = c\mathbf{e}^i$, where $\mathbf{e}^i_j =
\delta_{ij}$. It is obvious that $\Lambda(c \mathbf{e}^i) = \{0\}$
is a point. Then \eqref{eq:diag_action} becomes
\[
\Lambda(\v - c\mathbf{e}^i) \cong \Lambda(\v - c\mathbf{e}^i) \times
\Lambda(c\mathbf{e}^i) \stackrel{p_1}{\longleftarrow} \tilde
\Lambda(\v - c\mathbf{e}^i,c\mathbf{e}^i)
\stackrel{p_2}{\longrightarrow} \Lambda(\v).
\]

Assume $\Lambda(\v)_{i,c} \ne \emptyset$.  Then
\[
p_1^{-1}(\Lambda(\v - c\mathbf{e}^i)_{i,0}) =
p_2^{-1}(\Lambda(\v)_{i,c}).
\]
Let
\[
\tilde \Lambda(\v - c\mathbf{e}^i,c\mathbf{e}^i)_{i,0} =
p_1^{-1}(\Lambda(\v - c\mathbf{e}^i)_{i,0}) =
p_2^{-1}(\Lambda(\v)_{i,c}).
\]
We then have the following diagram.
\begin{equation}
\label{eq:crystal-action} \Lambda(\v - c\mathbf{e}^i)_{i,0}
\stackrel{p_1}{\longleftarrow} \tilde \Lambda(\v - c\mathbf{e}^i,
c\mathbf{e}^i)_{i,0} \stackrel{p_2}{\longrightarrow}
\Lambda(\v)_{i,c}
\end{equation}

\begin{lem}[\cite{KS97}]
\begin{enumerate}
\item The map $p_2 : \tilde \Lambda(\v - c\mathbf{e}^i,
c\mathbf{e}^i)_{i,0} \to \Lambda(\v)_{i,c}$ is a principal fiber
bundle with fiber $GL(\C^c) \times \prod_{i \in I} GL(V^{\v'}_i)$.

\item The map $p_1 : \tilde \Lambda(\v - c\mathbf{e}^i,
c\mathbf{e}^i)_{i,0} \to \Lambda(\v - c\mathbf{e}^i)_{i,0}$ is
smooth with fiber a connected rational variety.
\end{enumerate}
\end{lem}

\begin{cor}
\label{cor:irrcomp-isom} Suppose $\Lambda(\v)_{i,c} \ne \emptyset$.
Then there is a 1-1 correspondence between the set of irreducible
components of $\Lambda(\v - c\mathbf{e}^i)_{i,0}$ and the set of
irreducible components of $\Lambda(\v)_{i,c}$.
\end{cor}

Let $B(\v,\infty)$ denote the set of irreducible components of
$\Lambda(\v)$ and let $B_g(\infty) = \bigsqcup_\v B(\v,\infty)$. For
$X \in B(\v,\infty)$, let $\varepsilon_i(X) = \varepsilon_i(x)$ for
a generic point $x \in X$.  Then for $c \in \Z_{\ge 0}$ define
\[
B(\v,\infty)_{i,c} = \{X \in B(\v,\infty)\ |\ \varepsilon_i(X) =
c\}.
\]
Then by Corollary~\ref{cor:irrcomp-isom}, $B(\v -
c\mathbf{e}^i,\infty)_{i,0} \cong B(\v, \infty)_{i,c}$.

Suppose that ${\bar X} \in B(\v - c\mathbf{e}^i,\infty)_{i,0}$
corresponds to $X \in B(\v,\infty)_{i,c}$ by the above isomorphism.
Then we define maps
\begin{gather*}
\kf_i^c : B(\v - c\mathbf{e}^i,\w)_{i,0} \to B(\v,\w)_{i,c},\quad
\kf_i^c({\bar X})
= X, \\
\ke_i^c : B(\v,\w)_{i,c} \to B(\v - c\mathbf{e}^i,\w)_{i,0},\quad
\ke_i^c(X) = {\bar X}.
\end{gather*}
We then define the maps
\[
\ke_i, \kf_i : B_g(\infty) \to B_g(\infty) \sqcup \{0\}
\]
by
\begin{gather}
\ke_i : B(\v,\infty)_{i,c} \stackrel{\ke_i^c}{\longrightarrow} B(\v
- c\mathbf{e}^i, \infty)_{i,0}
\stackrel{\kf_i^{c-1}}{\longrightarrow} B(\v -
\mathbf{e}^i, \infty)_{i,c-1}, \label{ke-def} \\
\kf_i : B(\v,\infty)_{i,c} \stackrel{\ke_i^c}{\longrightarrow} B(\v
- c\mathbf{e}^i, \infty)_{i,0}
\stackrel{\kf_i^{c+1}}{\longrightarrow} B(\v + \mathbf{e}^i,
\infty)_{i,c+1}. \label{kf-def}
\end{gather}
We set $\ke_i(X)=0$ for $X \in B(\v,\infty)_{i,0}$.  Note that
$B(\v,\infty)_{i,c} \ne \emptyset$ implies $B(\v,\infty)_{i,c+1} \ne
\emptyset$ and thus $\kf_i(X)$ is never zero for $X \in
B(\v,\infty)_{i,c}$. We also define
\begin{gather*}
\wt : B_g(\infty) \to P,\quad \wt(X) = -\sum_{i\in I} \mathbf{v}_i
\alpha_i \text{
  for } X \in B(\v,\infty), \\
\varphi_i(X) = \varepsilon_i(X) + \left< h_i, \wt(X) \right>.
\end{gather*}
Here the $\alpha_i$ are the simple roots of $\g$.

\begin{prop}[\cite{KS97}]
$B_g(\infty)$ is a crystal and is isomorphic to the crystal
$B(\infty)$ of $U_q^-(\g)$.
\end{prop}


\section{Enumeration of components}
\label{sec:enumerate-comps}

In much the same way that Lusztig's quiver varieties yield a
geometric construction of the crystal graphs of
$U_q^-(\mathfrak{g})$, Nakajima's quiver varieties provide one with
a geometric realization of the crystal graphs of the irreducible
integrable highest-weight representations of $U_q(\mathfrak{g})$
(see \cite{S02}). In \cite{Sav03b}, an enumeration of the
irreducible components of Nakajima's quiver varieties for finite and
affine types $A$ and $D$ were given in terms of Young tableaux,
Young walls and new objects called Young pyramids.  In this section,
we describe a natural enumeration of the irreducible components of
Lusztig's quiver varieties in finite and affine type $A$.


\subsection{Type $A_n$}
\label{sec:a-comps}

In this subsection, let $\g = \mathfrak{sl}_{n+1}$ be the simple Lie
algebra of type $A_n$.  A key step in the enumeration of the
irreducible components of $\Lambda_\V$ is the following.

\begin{prop}[\cite{L91}]
\label{prop:irrcomp-finite} For \g\ a symmetric Lie algebra of
finite type (in particular, for $\g = \mathfrak{sl}_{n+1}$), the
irreducible components of $\Lambda_\V$ are the closures of the
conormal bundles of the various $G_\V$-orbits in $\mathbf{E}_{\V,
\Omega}$.
\end{prop}

By Proposition~\ref{prop:irrcomp-finite}, it suffices to enumerate
the $G_\V$-orbits in $\mathbf{E}_{\V,\Omega}$.  But these are simply
the isomorphism classes of the quiver $(I,\Omega)$.  By Gabriel's
Theorem, these are in one-to-one correspondence with the positive
roots of $\g$.  We can describe the isomorphism classes explicitly
as follows.

Let $I=\{1,\dots,n\}$ be the set of vertices of the Dynkin graph of
$\g$ with the set of oriented edges given by
\begin{gather*}
H=\{h_{i,j} \ |\ i,j \in I,\ |i-j|=1\}.
\end{gather*}
For two adjacent vertices $i$ and $j$, $h_{i,j}$ is the oriented
edge from vertex $i$ to vertex $j$.  Thus $\out(h_{i,j}) = i$ and
$\inc(h_{i,j})=j$. We define the involution $\bar{\ } : H \to H$ as
the function that interchanges $h_{i,j}$ and $h_{j,i}$. Let $\Omega
= \{h_{i,i-1}\ |\ 2 \le i \le n \}$.

For two integers $k,l$ such that $1 \le k \le l \le n$, define
$\V(k,l) \in \mathcal{V}$ to be the vector space with basis $\{ e_r\
|\ k \le r \le l\}$.  We require that $e_r$ has degree $r \in I$.
Let $x(k, l) \in \mathbf{E}_{\V(k,l), \Omega}$ be defined by $x(k,l)
: e_r \mapsto e_{r-1}$ for $k \le r \le l$, where $e_{k-1} = 0$. It
is clear that $(\V(k,l), x(k,l))$ is an indecomposable
representation of our quiver (i.e. element of
$\mathbf{E}_{\V,\Omega}$). Conversely, any indecomposable
finite-dimensional representation $(\V,x)$ of our quiver is
isomorphic to some $(\V(k,l),x(k,l))$. The correspondence guaranteed
by Gabriel's theorem is given by
\[
(\V(k,l),x(k,l)) \leftrightarrow \sum_{i=k}^{l} \alpha_i.
\]

Let $Z$ be the set of all pairs $(k,l)$ of integers such that $1 \le
k \le l \le n$ and let $\tilde Z$ be the set of all functions $Z \to
\N$ with finite support.

It is easy to see that for $\V \in \mathcal{V}$, the set of
$G_\V$-orbits in $\mathbf{E}_{\V, \Omega}$ is naturally indexed by
the subset $\tilde Z_\V$ of $\tilde Z$ consisting of those $\gamma
\in \tilde Z$ such that
\[
\sum_{(k,l) \in Z\, :\, k \le i \le l} \gamma(k,l) = \dim \V_i
\]
for all $i \in I$.  Corresponding to a given $\gamma$ is the orbit
consisting of all representations isomorphic to a sum of the
indecomposable representations $(V(k,l),x(k,l))$, each occurring
with multiplicity $\gamma(k,l)$. Denote by $\mathcal{O}_\gamma$ the
$G_\V$-orbit corresponding to $\gamma \in \tilde Z_\V$.  Let
$\mathcal{C}_f$ be the conormal bundle to $\mathcal{O}_\gamma$ and
let $\bar{\mathcal{C}}_\gamma$ be its closure. We then have the
following proposition.

\begin{prop}
The map $\gamma \mapsto \bar{\mathcal{C}}_\gamma$ is a one-to-one
correspondence between the set $\tilde Z_\V$ and the set of
irreducible components of $\Lambda_\V$.
\end{prop}
\begin{proof}
This follows immediately from Proposition~\ref{prop:irrcomp-finite}.
\end{proof}


\subsection{Type $A_n^{(1)}$}
\label{sec:ahat-comps}

In this subsection, we let $\g = \widehat{\mathfrak{sl}}_{n+1}$ be
the affine Lie algebra of type $A_n^{(1)}$.  Let $I=\Z/(n+1)\Z$ be
the set of vertices of the Dynkin graph of $\g$ with the set of
oriented edges given by
\begin{gather*}
H=\{ h_{i,j} \ |\ i,j \in I,\ i-j \equiv \pm 1 \mod n+1\}.
\end{gather*}

For two adjacent vertices $i$ and $j$, $h_{i,j}$ is the oriented
edge from vertex $i$ to vertex $j$.  Thus $\out(h_{i,j}) = i$ and
$\inc(h_{i,j})=j$. We define the involution $\bar{\ } : H \to H$ as
the function that interchanges $h_{i,j}$ and $h_{j,i}$. Let $\Omega
= \{h_{i,i-1}\ |\ i \in I \}$.

For two integers $k \le l$, define $\mathbf{V}(k,l) \in \mathcal{V}$
to be the vector space with basis $\{ e_r\ |\ k \le r \le l\}$.  We
require that $e_r$ has degree $r \mod n+1$. Let $x(k, l) \in
\mathbf{E}_{\mathbf{V}(k,l), \Omega}$ be defined by $x(k,l) : e_r
\mapsto e_{r-1}$ for $k \le r \le l$, where $e_{k-1} = 0$. It is
clear that $(\mathbf{V}(k,l), x(k,l))$ is an indecomposable
representation of our quiver and that $x(k,l)$ is nilpotent. Also,
the isomorphism class of this representation does not change when
$k$ and $l$ are simultaneously translated by a multiple of $n+1$.
Conversely, any indecomposable finite-dimensional representation
$(\mathbf{V},x)$ of our quiver, with $x$ nilpotent, is isomorphic to
some $(\mathbf{V}(k,l),x(k,l))$ where $k$ and $l$ are uniquely
defined up to a simultaneous translation by a multiple of $n+1$.

Let $Z$ be the set of all pairs $(k \le l)$ of integers defined up
to simultaneous translation by a multiple of $n+1$ and let $\tilde
Z$ be the set of all functions $Z \to \N$ with finite support.

It is easy to see that for $\mathbf{V} \in \mathcal{V}$, the set of
$G_\mathbf{V}$-orbits on the set of nilpotent elements in
$\mathbf{E}_{\mathbf{V}, \Omega}$ is naturally indexed by the subset
$\tilde Z_\mathbf{V}$ of $\tilde Z$ consisting of those $\gamma \in
\tilde Z$ such that
\[
\sum_{k \le l} \gamma(k,l) \cdot \#\{r\ |\ k \le r \le l,\ r \equiv
i \mod n+1\} = \dim \mathbf{V}_i
\]
for all $i \in I$.  Here the sum is taken over all $k \le l$ up to
simultaneous translation by a multiple of $n+1$.  Corresponding to a
given $\gamma$ is the orbit consisting of all representations
isomorphic to a sum of the indecomposable representations
$(V(k,l),x(k,l))$, each occurring with multiplicity $\gamma(k,l)$.
Denote by $\mathcal{O}_\gamma$ the $G_\mathbf{V}$-orbit
corresponding to $\gamma \in \tilde Z_\mathbf{V}$.

We say that $\gamma \in \tilde Z_\mathbf{V}$ is \emph{aperiodic} if
for any $k \le l$, not all $\gamma(k,l)$, $\gamma(k+1, l+1)$, \dots,
$\gamma(k+n, l+n)$ are greater than zero.  For any $\gamma \in
\tilde Z_\mathbf{V}$, let $\mathcal{C}_\gamma$ be the conormal
bundle of $\mathcal{O}_\gamma$ and let $\bar{\mathcal{C}}_\gamma$ be
its closure.

\begin{prop}[{\cite[15.5]{L91}}]
Let $\gamma \in \tilde Z_{\mathbf{V}}$.  The following two
conditions are equivalent.
\begin{enumerate}
\item $\mathcal{C}_\gamma$ consists entirely of
  nilpotent elements.
\item $\gamma$ is aperiodic.
\end{enumerate}
\end{prop}

\begin{prop}[{\cite[15.6]{L91}}]
The map $\gamma \to \bar{\mathcal{C}}_\gamma$ is a 1-1
correspondence between the set of aperiodic elements in $\tilde
Z_\mathbf{V}$ and the set of irreducible components of
$\Lambda_\mathbf{V}$.
\end{prop}


\section{Combinatorial crystal graphs}
\label{sec:comb-crystals}

In Section~\ref{sec:enumerate-comps}, we enumerated the irreducible
components of Lusztig's quiver varieties for finite and affine type
$A$. As noted in Section~\ref{sec:ca_qv}, there is a geometrically
defined crystal structure on the set of these irreducible
components. We now undertake the task of giving an explicit
description of this action in terms of the combinatorial data
enumerating irreducible components.  The arguments are similar to
those used in \cite{Sav03b}.

\subsection{Type $A_n$}
\label{sec:comb-action-a}

In this subsection, we consider the case where $\g =
\mathfrak{sl}_{n+1}$.  For $i \in I = \{1,\dots,n\}$, we define the
\emph{$i$-signature} of $\gamma \in \tilde Z$ as follows.  Consider
the ordering on pairs $(k,l)$, $1 \le k \le l \le n$, given by
\[
(k,l) < (k',l') \text{ if } k < k' \text{ or if } k=k' \text{ and }
l > l'.
\]
This is the lexicographic ordering where we reverse the usual order
on $\Z$ in the second factor.  Now, write all the pairs $(k,l)$ in
the domain of $\gamma$ from left to right in increasing order. Below
pairs of the form $(k,i-1)$, write $\gamma(k,i-1)$ $+$'s and below
pairs of the form $(k,i)$, write $\gamma(k,i)$ $-$'s.  Now, from
this list of $+$'s and $-$'s, cancel all $(+,-)$ pairs.  That is,
delete a $+$ and $-$ if they are adjacent in this list, with the $+$
to the left, and continue this process until so such configuration
exists. What is left is a (possibly empty) sequence of $-$'s
followed by a (possibly empty) sequence of $+$'s. We call this the
$i$-signature of $\gamma$.

Let $(k',i-1)$ be the pair corresponding to the leftmost $+$ in the
$i$-signature of $\gamma$ (that is, the pair under which this $+$
was placed in the above procedure).  If there is no $+$ in the
$i$-signature of $\gamma$, let $k'=i$.  Then define $\gamma^{i,+}$
by
\begin{align}
\gamma^{i,+}(k',i) &= \gamma(k',i) + 1, \\
\gamma^{i,+}(k',i-1) &= \gamma(k',i-1) - 1, \label{eq:+ignore}\\
\gamma^{i,+}(k,l) &= \gamma(k,l), \text{ for } k \ne k' \text{ or }
l \ne i,i+1.
\end{align}
In the case that $k'=i$, we ignore equation~\eqref{eq:+ignore}.

If there is at least one $-$ in the $i$-signature of $\gamma$, let
$(k'',i)$ be the pair corresponding to the rightmost $-$.  Then
define $\gamma^{i,-}$ by
\begin{align}
\gamma^{i,+}(k'',i) &= \gamma(k'',i) - 1, \\
\gamma^{i,+}(k'',i-1) &= \gamma(k'',i-1) + 1, \label{eq:-ignore}\\
\gamma^{i,+}(k,l) &= \gamma(k,l), \text{ for } k \ne k'' \text{ or }
l \ne i,i-1.
\end{align}
In the case that $k'=i$, we ignore equation~\eqref{eq:-ignore}.

\begin{theo}
\label{theo:a-comb-action}
For $\gamma \in \tilde Z$, let $X_\gamma
\in B_g(\infty)$ be the element of the crystal corresponding to the
irreducible component $\overline{\mathcal{C}}_\gamma$.  Then
\begin{align}
\varepsilon_i(X_\gamma) &= \# \{\text{$-$'s in the $i$-signature of
$\gamma$}\},
\label{eq:a-thm-epsilon} \\
\wt(X_\gamma) &= -\sum_{k \le i \le l}
f(k,l) \alpha_i, \label{eq:a-thm-wt} \\
\varphi_i(X_\gamma) &= \varepsilon_i(X_\gamma) +
\left<h_i,\wt(X_\gamma)\right>, \label{eq:a-thm-phi} \\
\ke_i(X_\gamma) &= \begin{cases} X_{\gamma^{i,-}} & \text{if }
\varepsilon_i(X_\gamma) > 0, \label{eq:a-thm-e} \\
0 & \text{if } \varepsilon(X_\gamma) = 0 \end{cases}, \\
\kf_i(X_\gamma) &= X_{\gamma^{i,+}}. \label{eq:a-thm-f}
\end{align}
\end{theo}
The remainder of this section is devoted to the proof of
Theorem~\ref{theo:a-comb-action}.

Consider the irreducible component $\overline{\mathcal{C}}_\gamma$
of $\Lambda_\V$.  This is the closure of the conormal bundle to the
orbit $\mathcal{O}_\gamma$.  Recall that for a point $x \in
\mathcal{O}_\gamma$, $\V$ decomposes into sums of copies of
$\V(k,l)$ for various $1 \le k \le l \le n$ corresponding to the
decomposition of $x$ into irreducible quiver representations
$x(k,l)$. Each copy of the form $\V(k,i)$ corresponds to a $-$ in
the construction of the $i$-signature of $\gamma$ (before the
$(+,-)$ cancelation was performed).

\begin{lem}
\label{lem:a-gen-point} For a generic point of the irreducible
component $\overline{\mathcal{C}}_\gamma$, $\gamma \in \tilde Z_\V$,
the image of $\oplus_{h\, :\, \inc(h)=i} x_h$ is spanned by the
degree $i$ vectors of those $\V(k,l)$, $k \le i < l$, appearing in
the decomposition of $\V$ and the degree $i$ vectors of those
$\V(k,i)$ appearing is this decomposition which correspond to $-$'s
that were canceled in the formation of the $i$-signature of
$\gamma$.
\end{lem}

\begin{proof}
The proof is almost identical to the proofs of Lemma~6.1 and
Proposition~6.2 of \cite{Sav03b} and thus is omitted.
\end{proof}

\begin{proof}[Proof of Theorem~\ref{theo:a-comb-action}]
It follows from Lemma~\ref{lem:a-gen-point} that
$\varepsilon_i(X_\gamma)$ is the number of $-$'s in the
$i$-signature of $\gamma$ and that
\[
\ke_i^{\varepsilon_i(X_\gamma)}(X_\gamma) = X_{\gamma'}
\]
where
\[
\gamma' = ( \cdots (\gamma^{i,-})^{i,-} \cdots ),
\]
and the superscript $i,-$ appears $\varepsilon_i(X_\gamma)$ times.

Now, it is easy to see from the definition of the action of $\kf_i$
on $B_g(\infty)$ that $\kf_i(X_\gamma)$ is never zero.  Then for
$\gamma \in \tilde Z$ with $\varepsilon(X_\gamma)=0$, we have
$\ke_i^c \kf_i^c (X_\gamma) = X_\gamma$.  Thus,
\[
\kf_i^c(X_\gamma) = X_{\gamma''}
\]
where
\[
\gamma'' = ( \cdots (\gamma^{i,+})^{i,+} \cdots ),
\]
and the superscript $i,+$ appears $c$ times.  Then
equations~\eqref{eq:a-thm-e} and~\eqref{eq:a-thm-f} follow
immediately.  Finally, equations~\eqref{eq:a-thm-wt}
and~\eqref{eq:a-thm-phi} follow from the general properties of
crystals.
\end{proof}


\subsection{Type $A_n^{(1)}$}

We now consider the case $\g = \widehat{\mathfrak{sl}}_{n+1}$.  Fix
an $i$ such that $0 \le i \le n$.  For elements $(k,l) \in Z$ such
that $l \equiv i$ or $i-1 \mod n+1$, pick the unique representative
such that $l=i$ or $i-1$.  Then define an ordering on such elements
by
\[
(k,l) < (k',l') \text{ if } k < k' \text{ or if } k=k' \text{ and }
l
> l'.
\]
We then mimic the definitions of Section~\ref{sec:comb-action-a}.
Write all the pairs of the form $(k,i)$ or $(k,i-1)$ (picking a
representative of this form for all pairs possible) in the domain of
$\gamma$ from left to right in increasing order.  Below pairs of the
form $(k,i-1)$, write $\gamma(k,i-1)$ $+$'s and below pairs of the
form $(k,i)$, write $\gamma(k,i)$ $-$'s.  From this list of $+$'s
and $-$'s, cancel all $(+,-)$ pairs.  What is left is a (possibly
empty) sequence of $-$'s followed by a (possibly empty) sequence of
$+$'s. We call this the $i$-signature of $\gamma$.

Let $(k',i-1)$ be the pair corresponding to the leftmost $+$ in the
$i$-signature of $\gamma$.  If there is no $+$ in the $i$-signature
of $\gamma$, let $k'=i$.  Then define $\gamma^{i,+}$ by
\begin{align}
\gamma^{i,+}(k',i) &= \gamma(k',i) + 1, \\
\gamma^{i,+}(k',i-1) &= \gamma(k',i-1) - 1, \label{eq:ahat+ignore}\\
\gamma^{i,+}(k,l) &= \gamma(k,l) \text{ for all other } (k,l) \in Z.
\end{align}
In the case that $k'=i$, we ignore equation~\eqref{eq:ahat+ignore}.

If there is at least one $-$ in the $i$-signature of $\gamma$, let
$(k'',i)$ be the pair corresponding to the rightmost $-$.  Then
define $\gamma^{i,-}$ by
\begin{align}
\gamma^{i,+}(k'',i) &= \gamma(k'',i) - 1, \\
\gamma^{i,+}(k'',i-1) &= \gamma(k'',i-1) + 1, \label{eq:ahat-ignore}\\
\gamma^{i,+}(k,l) &= \gamma(k,l) \text{ for all other } (k,l) \in Z.
\end{align}
In the case that $k'=i$, we ignore equation~\eqref{eq:ahat-ignore}.

\begin{theo}
\label{theo:ahat-comb-action} For an aperiodic $\gamma \in \tilde
Z$, let $X_\gamma \in B_g(\infty)$ be the element of the crystal
corresponding to the irreducible component
$\overline{\mathcal{C}}_\gamma$.  Then
\begin{align}
\varepsilon_i(X_\gamma) &= \# \{\text{$-$'s in the $i$-signature of
$\gamma$}\},
\label{eq:ahat-thm-epsilon} \\
\wt(X_\gamma) &= - \sum_{(k,l) \in Z} \sum_{i\, :\, k \le i \le l}
f(k,l) \alpha_{\bar i}, \label{eq:ahat-thm-wt} \\
\varphi_i(X_\gamma) &= \varepsilon_i(X_\gamma) +
\left<h_i,\wt(X_\gamma)\right>, \label{eq:ahat-thm-phi} \\
\ke_i(X_\gamma) &= \begin{cases} X_{\gamma^{i,-}} & \text{if }
\varepsilon_i(X_\gamma) > 0, \label{eq:ahat-thm-e} \\
0 & \text{if } \varepsilon(X_\gamma) = 0 \end{cases}, \\
\kf_i(X_\gamma) &= X_{\gamma^{i,+}}. \label{eq:ahat-thm-f}
\end{align}
\end{theo}
In the first sum in Equation~\eqref{eq:ahat-thm-wt}, we chose one
representative $(k,l)$ for each element in $Z$.  In this same
equation, $\bar i$ denotes the unique integer in $\{0,1,\dots,n\}$
congruent to $i$ modulo $n+1$.

\begin{proof}
The proof of this Theorem is similar to the proof of
Theorem~\ref{theo:a-comb-action} (see also \cite[Lemma~8.3 and
Theorem~8.4]{Sav03b}) and thus will be omitted.
\end{proof}


\section{Comparison to Young tableaux and Young walls}
\label{sec:comparison}

In Section~\ref{sec:comb-crystals} we gave explicit formulas for the
geometrically defined crystal structure on irreducible components of
quiver varieties in terms of the combinatorial enumerations of these
components described in Section~\ref{sec:enumerate-comps}.  There
are other combinatorial descriptions of the crystal $B(\infty)$ of
$U_q^-(\g)$ given in terms of Young tableaux and Young walls (see
\cite{HL05,Lee06,Lee05}).  In this section we show how our
combinatorial description can be translated into the Young
tableaux/wall approach and thus provide an alternative geometric
proof of these realizations in finite and affine type $A$.


\subsection{Type $A_n$}

Consider the case $\g = \mathfrak{sl}_{n+1}$.   We first recall the
realization of the crystal of $U_q^-(\g)$ in terms of Young tableaux
(see \cite{HL05,Lee05}).  A semi-standard tableau is called
\emph{large} if it consists of $n$ non-empty rows, and if for $1 \le
i \le n$, the number of $i$-entries in the $i$th row is strictly
greater than the number of all boxes in the $(i+1)$st row. Let
$\mathcal{T}^L$ denote the set of all large tableaux.

Two tableaux $T_1,T_2 \in \mathcal{T}^L$ are called \emph{related},
and we write $T_1 \sim T_2$, if for $1 \le i < j \le n$, the number
of $j$-entries in the $i$th rows of $T_1$ and $T_2$ are equal.  This
is an equivalence relation and we define $\mathcal{T}(\infty) =
\mathcal{T}^L/\sim$.  In \cite{Lee05},  a crystal structure on
$\mathcal{T}(\infty)$ is defined and it is shown that this crystal
is isomorphic to the crystal $B(\infty)$ for $U_q^-(\g)$.  The
crystal structure on $\mathcal{T}(\infty)$ is inspired by the
tableaux realization of irreducible integrable highest-weight
representations of $\g$ (see \cite{HK,KN94}).

Define a map $\tau : \mathcal{T}^L \to \tilde Z$ by setting
$\tau(T)(k,l)$, $1 \le k \le l \le n$, to be the number of
$(l+1)$-entries in the $k$th row of $T$. It follows immediately that
$\tau(T_1) = \tau(T_2)$ if $T_1 \sim T_2$. Thus, $\tau$ descends to
a map from $\mathcal{T}(\infty)$ to $\tilde Z$ which we also denote
by $\tau$.

\begin{theo}
For $\gamma \in \tilde Z$, let $X_\gamma$ be the element of
$B_g(\infty)$ corresponding to the irreducible component
$\overline{\mathcal{C}}_\gamma$.  The map $\mathcal{T}(\infty) \to
B_g(\infty)$ given by $T \mapsto X_{\tau(T)}$ is a crystal
isomorphism.
\end{theo}

\begin{proof}
Since the $\sim$-equivalence classes are uniquely determined by the
number of $j$-entries in the $i$th row for $j > i$, the map $\tau$
is one-to-one.  Now, it is easy to see that if $T_1 \sim T_2$ and $1
\le i \le n$, then the $i$-signatures of $T_1$ and $T_2$ differ only
by the possible addition of $+$'s to the right side.  These
additional $+$'s come from $i$-entries in the $i$th row. The action
of $\kf_i$ will act at one of these entries if and only if, there
are no $+$'s in the $i$-signature of $\gamma =
\gamma(T_1)=\gamma(T_2)$. But then the action of $\kf_i$ on $\gamma$
is to increase $\gamma(i,i)$ by one.  This corresponds to changing
an $i$-entry in the $i$th row to an $(i+1)$-entry in $T_1$ or $T_2$,
which is precisely how $\kf_i$ acts on these tableaux.  In the case
that $\kf_i$ does not act at one of these additional $+$'s, the fact
that the action of $\kf_i$ commutes with $\tau$ follows easily from
the definition of the action in terms of the $i$-signature.

The map $\varepsilon_i$, $1 \le i \le n$, commutes with $\tau$ since
it is simply the number of $-$'s in the $i$-signature of an element
of the crystal. The fact that $\wt$ and $\varphi_i$, $1 \le i \le
n$, commute with $\tau$ is also straightforward.
\end{proof}


\subsection{Type $A_n^{(1)}$}

We now let $\g = \widehat{\mathfrak{sl}}_{n+1}$.  A description of
the crystal $\mathcal{B}(\infty)$ is given in \cite{Lee06} in terms
of combinatorial objects called Young walls.  The set of all Young
walls built on the ``ground state'' wall $Y_\infty$ is denoted by
$\mathcal{F}(\infty)$ and the subset consisting of ``proper'' Young
walls is denoted by $\mathcal{Y}(\infty)$.  We define a map $\rho :
\mathcal{F}(\infty) \to \tilde Z$ as follows.

For $Y \in \mathcal{F}(\infty)$, each block that has been added to
$Y_\infty$ to form $Y$ (that is, appears in $Y$ but not $Y_\infty$)
sits to the left of a unique block in the rightmost slice which has
been added to $Y_\infty$. For each block in the rightmost slice of
$Y$ which has been added to $Y_\infty$, the number of blocks to its
left is finite.  For $0 \le k \le n$, we define $\rho(Y)(k,l)$ to be
the number of $(-k)$-colored blocks (recall that the set of block
colors is $\Z/(n+1)\Z$) in the rightmost slice of $Y$ which have
been added to $Y_\infty$ and which have a $l-k$ blocks to their
left. Each element of $Z$ has a unique representative $(k,l)$ with
$0 \le k \le n$ and thus $\rho(Y) \in \tilde Z$.  If $Y$ is reduced,
then $\rho(Y)$ is aperiodic. Thus $\rho(Y)$ restricts to a map from
$\mathcal{Y}(\infty)$ to the set of aperiodic elements of $\tilde
Z$.  This restriction is one-to-one.  Now, for a crystal $B$, define
a new crystal $B'$ by the map $\sigma : B \to B'$ that interchanges
the indices $i$ and $-i$ (this is induced by a Dynkin diagram
automorphism).  That is, $B' = B$ as sets and $\sigma$ is the
identity operators on the level of sets.  The crystal structure on
$B'$ is given by
\begin{align*}
\kf_i(\sigma(b)) &= \sigma(\kf_{-i}(b)), \\
\ke_i(\sigma(b)) &= \sigma(\ke_{-i}(b)), \\
\wt_i(\sigma(b)) &= \wt_{-i}(b), \\
\varepsilon_i(\sigma(b)) &= \varepsilon_{-i}(b), \\
\varphi_i(\sigma(b)) &= \varphi_{-i}(b).
\end{align*}
Recall that $\wt_i(b) = \left<h_i,\wt(b)\right>$. Note that while
$B'$ is a crystal with the above definitions, $\sigma$ is \emph{not}
a morphism of crystals in general.  For example, it does not commute
with the action of the operators $\ke_i$ and $\kf_i$.

\begin{theo}
For aperiodic $\gamma \in \tilde Z$, let $X_\gamma$ be the element
of $B_g(\infty)$ corresponding to the irreducible component
$\overline{\mathcal{C}}_\gamma$.  The map $Y \mapsto X_{\rho(Y)}$ is
a crystal isomorphism $\mathcal{Y}(\infty) \cong
\sigma(B_g(\infty))$.
\end{theo}

\begin{proof}
It is a technical but straightforward calculation to show that the
$i$-signatures of $Y$ and $X_{\rho(Y)}$ agree. Then the fact that
the map $Y \mapsto X_{\rho(Y)}$ commutes with the crystal operators
follows easily from their definitions in terms of the
$i$-signatures.
\end{proof}

Note that we could have avoided the appearance of the map $\sigma$
by choosing the opposite orientation of our quiver.  However, we
have chosen the orientation to agree with \cite{FS03,Sav03b}.  We
also note that while the proof that $\mathcal{Y}(\infty) \cong
B(\infty)$ found in \cite{Lee06} involves the Kyoto path realization
of $B(\infty)$, the alternative geometric proof presented here is
direct in the sense that it avoids any reference to this path
realization.  One enumerates the irreducible components of the
quiver varieties by Young walls and computes the action of the
crystal operators on these components directly.


\bibliographystyle{abbrv}
\bibliography{biblist}

\end{document}